\documentstyle[amstex,amscd,12pt]{article}
\makeatletter
\date{}
\newtheorem{theorem}{Theorem}[section]

 \newtheorem{conj}[theorem]{Conjecture}
\newtheorem{proposition}[theorem]{Proposition}

\newtheorem{problem}[theorem]{Problem}
\newcommand{\edim}{{\rm e-dim}}
\newcommand{\z}{{\Bbb Z}}
\newcommand{\q}{{\Bbb Q}}

\newcommand{\s}{{\Bbb S}}
\newcommand {\h}{{\check H}}

 \newcommand{\invlim}{\raisebox{-1ex}{$\stackrel{\hbox{lim}}{\leftarrow}$}}

\newcommand{\lo}{\longrightarrow}
\newcommand{\sm}{\setminus}
\newcommand{\tor}{{\rm Tor}}
\begin{document}

\title{Universal acyclic  resolutions for arbitrary coefficient groups\\
}

\author{Michael  Levin}

\maketitle
\begin{abstract}
We prove that for every compactum $X$   and every integer $n \geq 2$
there  are  a compactum $Z$ of $\dim \leq n+1$ and a surjective $UV^{n-1}$-map
$r: Z \lo X$   such that
 for every  abelian group $G$   and every integer $k \geq 2$
  such that  $\dim_G  X \leq k \leq n$ we have      $\dim_G Z \leq k$ and
  $r$ is $G$-acyclic.
\bigskip
\\
{\bf Math. Subj. Class.:} 55M10, 54F45.
\end{abstract}
\begin{section}{Introduction}
This paper is devoted to proving the following theorem which was announced in \cite{l4}.

\begin{theorem}
\label{t3}
Let $X$ be a compactum.  Then for every integer   $n \geq 2$
  there  are  a compactum $Z$ of $\dim \leq n+1$ and a surjective $UV^{n-1}$-map
$r: Z \lo X$   having the  property  that
 for every  abelian group $G$   and every integer $k \geq 2$
  such that  $\dim_G  X \leq k \leq n$ we have  that    $\dim_G Z \leq k$ and
  $r$ is $G$-acyclic.
   \end{theorem}

  The cohomological dimension $\dim_G X$
of  $X$  with respect to an abelian group $G$ is   the least number $n$ such that $\h^{n+1}
(X,A;G)=0$ for every closed subset $A$ of $X$.
 A space is    $G$-acyclic if its       reduced $\rm {\check C}$ech   cohomology groups modulo $G$
  are trivial,   a map is $G$-acyclic if
   every  fiber is  $G$-acyclic.
     By the Vietoris-Begle  theorem a surjective  $G$-acyclic map   of compacta
cannot raise the cohomological dimension $\dim_G$.
  A compactum $X$ is  approximately $n$-connected  if any embedding of $X$ into an ANR
  has the $UV^n$-property,   i.e.  for every    neighborhood $ U$ of $X$ there is
  a smaller neighborhood  $X \subset V \subset U$  such that the inclusion $V \subset U$
   induces
  the zero homomorphism of the homotopy groups in $\dim \leq n$.
  An     approximately $n$-connected  compactum has trivial    reduced $\rm {\check C}$ech
 cohomology groups     in $\dim \leq n$   with respect to any group $G$.
A map is called a    $UV^n$-map if  every fiber is      approximately $n$-connected.  \\

   Theorem \ref{t3} generalizes the following results of     \cite{l1, l3}.

        \begin{theorem}  {\rm (\cite{l1})}
\label{l1}
 Let $G$ be an abelian group and
 let $X$ be  a compactum with $\dim_G X \leq n$, $ n \geq 2$.  Then there are a compactum
 $Z $  with
 $\dim_G Z \leq n$ and  $\dim Z \leq n+1$ and a $G$-acyclic map $r : Z \lo X$  from $Z$ onto $X$.
\end{theorem}

   \begin{theorem}   {\rm(\cite{l3})}
 \label{l3}
 Let $X$ be a compactum  with $\dim_\z X \leq n \geq 2$.  Then there exist
 a compactum $Z$ with $\dim Z \leq n$  and a cell-like map $r : Z \lo X$ from
 $Z$ onto $X$ such for every integer $k \geq 2$ and every group $G$  such that
 $\dim_G X \leq k$ we have $\dim_G Z \leq k$.
 \end{theorem}

    Theorem \ref{l1} obviously follows from Theorem \ref{t3}.    Theorem \ref{l3} can be derived
  from      Theorem \ref{t3}   as follows.
      Recall that a   compactum     is cell-like if any map  from the compactum  to a CW-complex   is
   null-homotopic.
   A  map is cell-like if its fibers are cell-like.
    Let $X$ be of $\dim_\z  < \infty$  and let
  $r : Z \lo X$  satisfy the conclusions of  Theorem \ref{t3}  for $n= \dim_\z X +1$.
  Then $\dim_\z Z \leq \dim_\z X \leq n-1$ and because $Z$ is finite dimensional we have
  $\dim Z =\dim_\z Z \leq n -1$.     Since  $r$ is $UV^{n-1}$     and $\dim Z \leq n-1$
  we get that $r$ is cell-like.  Let for a group $G$, $\dim_G X \leq k \geq 2$.
   If $k \leq n$ then  by    Theorem \ref{t3},  $\dim_G  Z\leq k$ and if   $k > n$
   then $ \dim_G  Z\leq k$    since $\dim Z \leq n-1$. Thus
       Theorem \ref{t3}  implies      Theorem \ref{l3}.   \\

       It was observed in \cite{l3} that the restriction $k \geq 2$  in Theorem \ref{l3}
       cannot be omitted.  Therefore  Theorem \ref{t3}
       does not hold for $k=1$.

Let us discuss possible generalizations of Theorem \ref{t3}.
  One is tempted to reduce the dimension of $Z$ to $n$.
It is  partially  justified
by

\begin{theorem}
{\rm (\cite{l4})}
\label{t1}
 Let $X$ be a compactum.  Then for every integer   $n \geq 2$
  there  are  a compactum $Z$ of $\dim \leq n$ and a surjective $UV^{n-1}$-map
$r: Z \lo X$   having the  property  that
 for every finitely generated abelian group $G$   and every integer $k \geq 2$
  such that  $\dim_G  X \leq k \leq n$ we have  that    $\dim_G Z \leq k$ and
  $r$ is $G$-acyclic.
  \end{theorem}

  However Theorem \ref{t1} does not hold for arbitrary groups $G$.  Indeed,  one can show
  that      a  $\q$-acyclic  $UV^1$-map  from
a  compactum  of $\dim \leq 2$ must be  $\z$-acyclic (even cell-like).
  Then a compactum $X$ with $\dim_\z =3$ and $\dim_\q =2$ cannot be the image
  of a compuctum of $\dim \leq 2$ under     a  $\q$-acyclic  $UV^1$-map.

  The situation becomes more complicated if we drop in Theorem
  \ref{t3}  the requirement that $r$ is   $UV^{n-1}$ and consider

  \begin{problem}
  \label{pr}
  Given a compactum $X$,   an integer $n\geq 2$ and
  a collection of abelian groups ${\cal G}$  such that  $\dim_G X  \leq n$
  for every   $G \in \cal G$
  do there exist a compactum $Z$ of $\dim \leq n$ and a $\cal G$-acyclic  surjective map
  $r : Z \lo X$ such that $\dim_G Z \leq  \max\{\dim_G X , 2\}$
  for every $G \in \cal G$?  (The $\cal G$-acyclicity means the  $G$-acyclicity
   for every   $G \in \cal G$.)
  \end{problem}

In general   the answer to  Problem \ref{pr} is negative \cite{ky2}.
 Indeed,
let $X$ be a compactum    with $\dim_\z =3$, $\dim_\q =2$ and  $\dim_{\z_p}=2$ for every prime $p$
and let $G=\q \oplus (\q/\z)$.  Clearly  $\dim_G X =2$ and
 the $G$-acyclicity implies both the $\q$ and $(\q/\z)$-acyclicities.
Then it follows from the Bockstein sequence generated by

$0 \lo \z   \lo \q \lo \q /\z \lo 0$
 \\
 that    the $G$-acyclicity implies   the  $\z$-acyclicity
and therefore there is no $ G$-acyclic resolution for $X$ from a compactum of $\dim \leq  2$.

The situation described in the example can be interpreted in terms of Bockstein Theory.
Let $\cal G$ be a collection of abelian groups.
Denote by  $\sigma({\cal G}) $    the union   of the Bockstein basises
$ \sigma(G)$  of all  $G \in   { \cal G} $.
Based on the Bockstein inequalities define the closure $\overline{\sigma({\cal G})}$
of $\sigma({\cal G}) $  as
a collection of abelian groups containing ${\sigma({\cal G})}$    and possibly     some
 additional groups
determined by:

$\z_p   \in \overline{\sigma({\cal G})}$   if $\z_{p^\infty} \in     \overline{\sigma({\cal G})}$;

$\z_{p^\infty}      \in \overline{\sigma({\cal G})}$   if $\z_{p} \in     \overline{\sigma({\cal G})}$;

$\q \in     \overline{\sigma({\cal G})}$  if     $\z_{(p)} \in     \overline{\sigma({\cal G})}$;

$\z_{(p)}        \in \overline{\sigma({\cal G})}$
if $\q $ and $\z_{p^\infty} \in     \overline{\sigma({\cal G})}$.
\\
One can show that for compact metric spaces the  ${\cal G}$-acyclicity implies
 the $  \overline{\sigma({\cal G})}$-acyclicity.
This motivates the following

\begin{conj}
Problem \ref{pr}  can be  answered positively     if    $\dim_E X \leq n$ for every
$ E   \in \overline{\sigma({\cal G})}$.
\end{conj}

The  key open case of this conjecture seems  to be  constructing a $\q$-acyclic resolution
$r : Z \lo X$  for a  compactum $X$ with $\dim_\q \leq n$, $n\geq 2$  from a compactum $Z$
of $\dim \leq n$.

 \end{section}
 \begin{section}{Preliminaries}
  All  groups are assumed to be
   abelian  and functions between groups are homomorphisms.
   $\cal P$ stands for the set of  primes.  For a non-empty subset $\cal A$ of $ \cal P$ let
$S({\cal A})=\{ p_{1}^{n_1}p_{2}^{n_2}...p_{k}^{n_k} : p_i \in {\cal A}, n_i \geq 0\}$
 be  the set of  positive  integers with prime factors from $\cal A$ and for the empty set define
 $S(\emptyset)=\{ 1 \}$.
Let   $G$ be a group and  $g \in G$.
We say that $g$  is $\cal A$-torsion  if there is $n\in S(\cal A)$
such that $ng=0$  and  $g$  is $\cal A$-divisible if for every $n \in S(\cal A)$ there is $h\in G$
such that $nh=g$.
$\tor_{\cal A} G$ is the subgroup of the $\cal A$-torsion elements
of $G$.   $G$ is  $\cal A$-torsion  if $G=\tor_{\cal A}   G$, $G$ is $\cal A$-torsion free
if $\tor_{\cal A} G= 0$
and $G$ is $\cal A$-divisible if every element of $G$ is $\cal A$-divisible.

$G$ is $\cal A$-local if $G$ is $({\cal P} \sm {\cal A})$-divisible
and       $({\cal P} \sm {\cal A})$-torsion free.    The  $\cal A$-localization of  $G$ is
the homomorphism  $G \lo G \otimes \z_{(\cal A)}$  defined by $g \lo g \otimes 1$ where
$\z_{(\cal A)}= \{ n/m : n \in \z, m\in S({\cal P}\sm {\cal A})\}$.     $G$ is  $\cal A$-local
if and only if the $\cal A$-localization of  $G$  is an isomorphism.
A simply connected
CW-complex is $\cal A$-local if its homotopy groups are $\cal  A$-local.
A map between
two simply connected
 CW-complexes is  an $\cal A$-localization  if the induced homomorphisms of
the homotopy  and the (reduced integral)  homology groups are $\cal    A$-localizations.      \\

The extensional dimension  of a compactum
 $X$      is said not to exceed      a CW-complex $K$
 $K$, written $\edim X \leq K$,  if for every closed subset $A$
 of $X$ and every map $f : A \lo K$ there is an extension of $f$ over $X$.
 It is well-known that $\dim X \leq n$ is equivalent to $\edim X \leq \s^n$
 and $\dim_G X \leq  n$ is equivalent to $\edim X \leq K(G,n)$ where
 $K(G,n)$ is an Eilenberg-Mac Lane complex of type $(G,n)$. \\

     A  map between CW-complexes  is said to be  combinatorial if  the preimage of
   every subcomplex
   of the range is a subcomplex of the domain.
   Let $M$ be  a simplicial complex and let  $M^{[k]}$        be
    the $k$-skeleton of $M$ (=the union of all simplexes of $M$ of $\dim \leq k$).
By
 a resolution $EW(M,k)$   of $M$   we mean a CW-complex $EW(M,k)$ and
 a combinatorial map
 $\omega : EW(M,k) \lo  M$ such that $\omega$ is 1-to-1 over $M^{[k]}$.
 Let $f : N \lo K$  be a map of a subcomplex $N$ of $M$ into a CW-complex $K$.
The resolution is said to be suitable for $f$   if
 the map  $f \circ\omega|_{\omega^{-1}(N)}$ extends
   to a map  $ f': EW(M,k) \lo K$.  We call $f'$  a resolving map for $f$.
     The resolution is said to be
   suitable  for  a compactum $X$
if for   every simplex $\Delta$ of $M$,
 $\edim X  \leq \omega^{-1}(\Delta)$.
     Note that if $\omega: EW(M,k) \lo M$ is a resolution suitable
 for $X$   then  for every map  $\phi :  X \lo  M$  there is  a map   $\psi : X \lo EW(M,k)$
 such that  for every simplex $\Delta$ of $M$,
  $(\omega \circ \psi)(\phi^{-1}(\Delta)) \subset \Delta$.
 We  call $\psi$ a combinatorial lifting of $\phi$.

Let $M$   be a finite simplicial complex and
 let       $f : N  \lo K$ be a cellular  map  from a subcomplex $N$
of $M$  to a CW-complex    $K$ such that
$M^{[k]}\subset N$.
A standard way of constructing
a  resolution  suitable for $f$ is described in \cite{l3}.
Such a  resolution  $\omega: EW(M,k) \lo M$  is called the standard resolution    of $M$
for $f$ and it has the following properties:

$EW(M,k)$ is $(k-1)$-connected if so are $M$ and $K$;

$\omega$ is a map onto and for every simplex $\Delta$ of $M$, $\omega^{-1}(\Delta)$ is  either
contractible  or homotopy equivalent to $K$;

for every subcomplex $T$ of $M$,  $\omega|_{\omega^{-1}(T)} : EW(T,k)= \omega^{-1}(T) \lo T$
is the standard resolution of $T$ for $f|_{N \cap T} :   N \cap T \lo K$.
     \\

     Let $G$ be a group,  let $\alpha : L \lo M$ be a surjective  combinatorial map of a CW-complex
$L$ and a finite simplicial complex $M$ and let $n$ be a positive integer
 such that  $\Tilde{H}_i(\alpha^{-1}(\Delta) ; G)=0$
for every $i < n$ and  every  simplex $\Delta $ of $M$.     One can show by induction on
the number of simplexes of $M$  using  the Mayer-Vietoris sequence and the Five Lemma
that $\alpha_* :   \Tilde{H}_i (L; G) \lo    \Tilde{H}_i (M; G)$ is an isomorphism
for $i < n$.   We will refer to this fact as the combinatorial Vietoris-Begle theorem.  \\

 \begin{proposition}
 \label{p-dc}
 {\rm (\cite{dc})}. Let $G$ be a group and $p \in \cal P$.
 The following conditions are equivalent:

 $G$ is $p$-divisible;

 ${\rm Ext}(\z_{p^\infty}, G)$ is $p$-divisible;

 $   {\rm Ext}(\z_{p^\infty}, G)=0$.
 \end{proposition}

\begin{proposition}
\label{p1}
Let $G$ be a group,  let  $2 \leq k \leq n $  be  integers and let  $\cal F\subset \cal P$ and
$p \in {\cal P }\sm \cal F$.      Let $M$ be an $(n-1)$-connected fine simplicial    complex
such that
 $H_n(M)$ is $\cal F$-torsion  and
let     $\omega : L =EW(M,k) \lo M$ be       the standard resolution of $M$ for
    a cellular map  $f : N \lo K(G,k)$  from   a subcomplex $N$  of $M$ containing $M^{[k]}$.
 Then
 $L$  is $(k-1)$-connected and for every $1 \leq i \leq n-1$

 (i) $\pi_i (L)$ and $\pi_n (L)/ \tor_{\cal F}\pi_n (L)$    are $p$-torsion if $G=\z_p$;

 (ii)      $\pi_i (L)$ and $\pi_n (L)/ \tor_{\cal F}\pi_n (L)$    are $p$-torsion and
  $\pi_k (L)$ is $p$-divisisible if $G=\z_{p^\infty}$

 (iii)       $\pi_i (L)$ and $\pi_n (L)/ \tor_{\cal F}\pi_n (L)$    are
 $q$-divisible  and   $\pi_i (L)$  is $q$-torsion free
 for every $q\in \cal P$, $ q \neq  p$  if     $G=\z_{(p)}$;

 (iv)       $\pi_i (L)$ and $\pi_n (L)/ \tor_{\cal F}\pi_n (L)$    are
 $q$-divisible   and     $\pi_i (L)$  is $q$-torsion free
 for every $q\in \cal P$   if     $G=\q$.

\end{proposition}
  {\bf Proof.}
  Recall that   $\omega$ is a combinatorial surjective map,
    for every  simplex $\Delta$ of $ M$,
  $\omega^{-1}(\Delta)$ is either contractible  or homotopy equivalent to  $K(G,k)$
  and
 $L$ is $(k-1)$-connected because  so are $M$ and        $K(G,k)$.
 Since $M$ is    $(n-1)$-connected and   $H_n(M)$ is $\cal F$-torsion we have   that
$ H_n (M;\q)=0$ and $ H_n( M; \z_q)=0$, $ H_n( M; \z_{(q)})=0$
  for
   $q \in {\cal P} \sm \cal F$  and $ H_n( M; \z_{q^\infty})=0$  for every $q \in \cal P$.\\

(i)
 By the generalized  Hurewicz theorem
  $\tilde H_*( K(\z_p, k))$
  is  $p$-torsion. Then    $ \tilde H_*( K(\z_p, k); \q)=0$.
  Hence  by
  the combinatorial Vietoris-Begle theorem       $ \tilde H_i( L; \q)=0$  for $i \leq n$
  and therefore      $ \tilde H_i( L)$   is torsion for  $i \leq n$.
  \\
  \\
  Let $q\in \cal P$ and $q \neq p$  and $i \leq n-1$.\\
  Note  $ \tilde H_*( K(\z_p, k); \z_{(q)})=0$ and hence
  by the combinatorial Vietoris-Begle theorem       $ \tilde H_i( L; \z_{(q)})=0$.
  Then  $ \tilde H_i( L) \otimes \z_{(q)}=0$. Thus   $ \tilde H_i( L)$ is torsion
  and  $q$-torsion free
and    hence $ \tilde H_i( L)$    is $p$-torsion.
\\
\\
Now let  $q\in {\cal P}\sm \cal F$
 and $q \neq p$.\\
       Recall  $ H_n(M; \z_{(q)})=0$. Then       using the previous argument
       we conclude that
 $  H_n( L)$ is $q$-torsion free  and hence  $  H_n( L)$ is $({\cal F}\cup \{p\})$-torsion. \\

  By  the  generalized  Hurewicz theorem
   $\pi_i (L)$  is $p$ torsion for $i \leq n-1$ and        $\pi_n(L)$         is $({\cal F}\cup \{p\})$-torsion.
  Thus    $\pi_n (L)/ \tor_{\cal F}\pi_n (L)$  is $p$-torsion and (i) follows.  \\

  (ii)     The argument  used in (i) applies to show that
         $\pi_i (L)$,  $i \leq n-1$  and $\pi_n (L)/ \tor_{\cal F}\pi_n (L)$    are $p$-torsion.
  Note that  $\pi_k(L)=   H_k(L)$.  We will show
    that      $H_k(L)$ is $p$-divisible  and  this will imply (ii).
     Observe that
    $H_k( K(\z_{p^\infty},k);\z_p)=  \z_{p^\infty} \otimes  \z_p=0$.
    Then since $H_k (M;\z_p)=0$ the combinatorial  Vietoris-Begle theorem implies that
    $H_k (L;\z_p)=0$.
    Thus  $H_k( L)\otimes \z_p =0$ and therefore    $H_k( L)$ is $p$-divisible. \\

    (iii)
Since   $Z_{(p)}$ is $p$-local     we have that
    $\tilde H_*(K(Z_{(p)}, k))$ is    $p$-local    and therefore
    $\tilde H_*(K(Z_{(p)}, k); \z_q)=\tilde H_*(K(Z_{(p)}, k); \z_{q^\infty})=0$
    for every     $q \in \cal P$ , $q\neq p$.
  \\
   \\   Let    $q \in \cal P$ , $q\neq p$.\\
 Recall that $ \tilde H_i(M; \z_{q^\infty})=0$ for $i \leq n$.
      Then by  the combinatorial  Vietoris-Begle theorem
    $\tilde H_i (L;\z_{q^\infty})=0$ for  $i \leq n$.
    Hence by virtue
   of the universal coefficient  theorem         $\tilde H_i   ( L)*\z_{q^\infty} =0$     and
   $\tilde H_i   ( L)\otimes \z_{q^\infty} =0$
   for $i \leq n-1$  and therefore    $\tilde H_i   ( L)$ is $q$-torsion free and
   $q$-divisible
   for $i \leq n-1$.
     \\
     \\
   Let $q \in \cal P$ , $q\neq p$ and $q \notin \cal F$.
   \\
   Recall that $  H_n(M; \z_{q})=0$.
   By  the combinatorial  Vietoris-Begle theorem
   $ H_n (L;\z_{q})=0$.
   Hence $H_n   ( L)\otimes \z_q =0$       and        therefore     $ H_n   ( L)$ is $q$-divisible.
   \\
   \\
 Let  $q \in \cal P$ , $q\neq p$ and $q \in \cal F$.
 \\
 Then
$ H_n(M; \z_{q^\infty})=0$.
   By  the combinatorial  Vietoris-Begle theorem $ H_n (L;\z_{q^\infty})=0$.
  Hence $ H_n   ( L)\otimes \z_{q^\infty} =0$       and        therefore
  $H_n   ( L)/\tor_{q}H_n   ( L) $ is $q$-divisible.
   \\
   \\
   Now using  Completion and Localization Theories \cite{bk} we will pass  to the homotopy groups
   of $L$.
   \\  \\
        Let $q \in \cal P$ and $q\neq p$.  \\
         Denote   ${\cal A}= {\cal P} \sm \{ q\}$ and   let
          $\alpha : L \lo L_\alpha$ be an  $\cal A$-localization
  of $L$.
 Recall that  $\tilde H_i   ( L)$ is $q$-torsion free and
   $q$-divisible   for $i\leq n-1$. Then
   $\alpha$ induces an isomorphism of the  groups   $\tilde{H}_*   ( L)$
  and $\tilde{H}_*   ( L_\alpha )$  in  $\dim \leq n-1$.
Note that  $H_n(L) \otimes \z_{(\cal A)}= H_n(L) \otimes
\z{[1/q]}=
    (H_n   ( L)/\tor_{q}H_n   ( L))     \otimes \z{[1/q]}$ and  since
$ H_n   ( L)/\tor_{q}H_n   ( L)$   is $q$-divisible we have  that
the $\cal A$-localization of $  H_n   ( L)$ is an epimorphism.
Then        by the Whitehead theorem $\alpha$ induces  an
isomorphism of $\pi_i(L)$ and $\pi_i(L_\alpha)$, $i\leq n-1$ and
an epimorphism of $\pi_n(L)$ and $\pi_n(L_\alpha)$. Therefore
$\pi_i(L),i\leq n-1$ is $\cal A$-local (that is, $q$-divisible and
$q$-torsion free) and   the $\cal A$-localization of $\pi_n(L)$ is
an epimorphism. The last property implies
 that $\pi_n ( L)/\tor_{q}\pi_n   ( L)$ is $q$-divisible.
     \\\\
   Let   $q \in \cal P$ , $q\neq p$ and    $q \notin \cal F$.\\
   Let $\beta  : L \lo  L_\beta$ be        a     $q$-completion of $L$.
   Then $\beta$ induces an isomorphism of $\tilde{H}_*(L;\z_q)$ and   $\tilde{H}_*(L_\beta ;\z_q)$
   and since  $H_n(L;\z_q)=0$  we get that    $H_n(L_\beta ;\z_q)=0$
   and therefore      $H_n(L_\beta)$ is $q$-divisible.   Now since      $\pi_i(L)$ is
   $q$-divisible and $q$-torsion free      we have that       ${\rm Hom}(\z_{q^\infty},   \pi_i(L))=0$,
 $i \leq n-1$
   and  by Proposition \ref{p-dc}  ${\rm Ext}(\z_{q^\infty},   \pi_i(L))=0$, $i \leq n-1$.
   Then the exact sequence   \\

 $0\lo    {\rm Ext}(\z_{q^\infty},   \pi_i(L)) \lo \pi_i (L_\beta) \lo
   {\rm Hom}(\z_{q^\infty},   \pi_{i-1}(L)) \lo 0$\\
   \\
   implies  that $L_\beta$ is $(n-1)$-connected and
   $ {\rm Ext}(\z_{q^\infty},   \pi_n(L)) = \pi_n (L_\beta)   $.  Thus $ \pi_n (L_\beta) =H_n(L_\beta)$
      and hence     $ {\rm Ext}(\z_{q^\infty},   \pi_n(L))$ is $q$-divisible.  Then
      by Proposition \ref{p-dc}  $ \pi_n(L)$  is $q$-divisible and
   (iii) is proved.     \\    \\

 (iv).  The proof is similar to the proof of (iii).
\hfill $\Box$
     \\\\
    Let $X$ be a compactum  and let $n$ be a positive integer.
   The Bockstein basis  of abelian groups is
 the following collection  of groups $\sigma = \{\q,  \z_p,   \z_{p^\infty}, \z_{(p)}    : p \in {\cal P} \}$.
 Define  the Bockstein basis
 $\sigma(X,n)$ of $X$ in dimensions $\leq n$ as
 $\sigma(X,n) =\{ E \in \sigma : \dim_E X \leq n \}$.   Following \cite{l1}
    denote:\\

   ${\cal T }(X,n)= \{ p \in{ \cal P} :  \z_p $ or $\z_{p^\infty} \in \sigma (X,n)\}$ ; \\

   ${\cal D }(X,n)=\emptyset$  if      $ \sigma (X,n)$   contains only torsion groups;

  ${\cal D }(X,n)=  \cal P$ if $\q \in \sigma (X,n)$ but none of $\z_{(p)}$, $p\in\cal P$
  belongs to     $ \sigma (X,n)$  and

  ${\cal D }(X,n)={\cal P} \sm  \{ p \in{ \cal P} :  \z_{(p)} \in \sigma (X,n) \}$       otherwise  ;\\

    ${\cal F}(X,n)={\cal D}(X,n) \sm {\cal T}(X,n) $. \\
\\
   Note  that
   for every group $G$ such that $\dim_G X \leq n$,
   $G$ is   ${\cal F}(X,n)$-torsion free.

    \begin{proposition}
\label{h3}
 Let $X$ be a compactum   such that      ${\cal D }(X,n) \neq \emptyset$.
 Then
   $\dim_H X \leq n$ for  every  group  $H$
   such that $H$ is ${\cal D}(X,n)$-divisible and ${\cal F}(X,n)$-torsion free.

\end{proposition}
{\bf Proof.}  Let $G=\oplus\{ E  : E \in \sigma(X,n) \}$.  Then
$\dim_G X \leq n$. One can easily verify that in the notations of Proposition 2.4 of \cite{l1},
$ {\cal D}(G)=  {\cal D}(X,n)$ and      $ {\cal F}(G)=  {\cal F}(X,n)$.
Then the result follows from      Proposition 2.4 of \cite{l1}.
 \hfill $\Box$    \\

 In the proof of Theorem \ref{t3}    we will also  use the   following facts.
  \begin{proposition}  {\rm (\cite{l3})}
 \label{p2}
 Let $K$ be a simply connected CW-complex such that  $K$ has only finitely many non-trivial
 homotopy groups.   Let $X$ be a compactum such that $\dim_{\pi_i (K)} X  \leq i$
 for  $i>1$.  Then $\edim X \leq K$.
 \end{proposition}
 Let $K'$ be a simplicial complex.  We say that   maps
  $h : K \lo K'$, $g : L \lo L'$,  $\alpha : L \lo K$ and $\alpha' : L' \lo K'$
  combinatorially commute   if   for every
 simplex  $\Delta$ of $K'$ we have that
 $(\alpha' \circ g)(( h\circ \alpha)^{-1}(\Delta)) \subset \Delta$.
 Recall that a  map $h' : K \lo L'$ is
 a combinatorial lifting of $h$  to $L'$  if
 for every  simplex  $\Delta$ of $K'$   we have that
 $(\alpha' \circ h')( h^{-1}(\Delta)) \subset \Delta$.

 For   a simplicial complex $K$ and $a \in K$,   $st(a)$  denotes the union of all the simplexes
 of   $K$ containing  $a$.

 \begin{proposition} {\rm (\cite{l3})}
 \label{p3}${}$

   (i)     Let a compactum $X$  be  represented as the inverse limit
 $X  ={\rm  \invlim }K_i$ of  finite simplicial complexes $K_i$
 with  bonding maps   $h_{j}^i  :  K_{j} \lo K_i$.  Fix $i$ and   let
 $\omega : EW(K_i,k) \lo K_i$ be a resolution of $K_i$ which is suitable for $X$.
 Then there is
  a sufficiently large $j$  such that $h_j^i$ admits a combinatorial lifting to $EW(K_i,k)$.

  (ii)       Let  $h : K \lo K'$, $h' : K \lo L'$ and $\alpha' : L' \lo K'$ be maps of
  a simplicial complex $K'$   and  CW-complexes $K$ and $L'$  such that
  $h$ and $\alpha'$ are combinatorial and $h'$ is a combinatorial lifting of $h$.
 Then  there is a cellular approximation of $h'$ which is  also a combinatorial lifting of $h$.

 (iii)    Let $K$ and $K'$ be simplicial complexes,
 let    maps    $h : K \lo K'$, $g : L \lo L'$,  $\alpha : L \lo K$ and $\alpha' : L' \lo K'$
 combinatorially commute and let $h$ be  combinatorial.  Then

 $ g( \alpha^{-1}(st( x)) ) \subset
   \alpha'{}^{-1}(st( h (x)) )$ and
   $h(st((\alpha (z)))  \subset        st((\alpha' \circ g)(z))$  \\
   for every   $x \in K$ and  $z\in L$;
 \end{proposition}

  \end{section}

  \begin{section}{Proof of Theorem \ref{t3}}
  Denote ${\cal D}={\cal D}(X, n)$   and      ${\cal F}={\cal F}(X, n)$.
     Represent $X$ as   the inverse limit $X  = \invlim (K_i,h_i)$ of finite simplicial complexes $K_i$
 with combinatorial bonding maps   $h_{i+1}  :  K_{i+1} \lo K_i$  and the projections
 $p_i : X \lo K_i$ such that for every simplex $\Delta$ of $K_i$,  diam$(p_i^{-1}(\Delta)) \leq 1/i$.
 Following A. Dranishnikov \cite{d0.5, d1} we construct  by induction finite  CW-complexes
  $L_i$   and maps      $g_{i+1}: L_{i+1} \lo L_i$,
 $\alpha_i : L_i \lo K_i$   such that  \\
\\
(a)  $L_i$ is $(n+1)$-dimensional  and
  obtained from $K^{[n+1]}_i$ by replacing some $(n+1)$-simplexes by $(n+1)$-cells
 attached  to the   boundary of the  replaced simplexes by a map of degree $\in  S({\cal F})$.
 Then $ \alpha_i$ is a projection of $L_i$ taking
 the new cells to the  original ones such that $\alpha_i$ is 1-to-1 over
 $K^{[n]}_i$.  We  define
 a simplicial structure on $L_i$ for which $\alpha_i$ is a  combinatorial map and refer
 to this simplicial structure while constructing resolutions of $L_i$.
 Note
 that for  ${\cal F}=\emptyset$ we don't replace simplexes of $K^{[n+1]}_i$ at all;
  \\
  \\
 (b)  the maps $h_{i+1}$, $g_{i+1}$,  $\alpha_{i+1}$  and $\alpha_i$   combinatorially commute.
 Recall that   this  means that
 for every simplex  $\Delta $ of $ K_{i}$,
 $(\alpha_i  \circ g_{i+1})((h_{i+1} \circ \alpha_{i+1})^{-1}(\Delta)) \subset \Delta$.  \\
  \\
  We will construct $L_i$ in such a way that   $Z=\invlim(L_i,g_i)$ will  admit
a map $r : Z \lo X$ such that  $Z$ and $r$ satisfy the conclusions of the theorem.

   Let     $E \in    \sigma$  be such that
 $\dim_{E}  X \leq k$, $2 \leq k  \leq n$ and let $f  :  N  \lo K(E ,k)$ be a  cellular map  from
 a subcomplex $N$ of $L_i$, $  L_i^{[k]}  \subset  N$.
Let  $\omega_L: EW(L_i, k) \lo L_i^{}$
be the standard resolution  of $L_i$ for $f$.   We are  going to construct
from  $\omega_L: EW(L_i, k) \lo L_i^{}$ a resolution $\omega : EW(K_i,k)\lo K_i$ of
 $K_i$ suitable for $X$.  If $\dim K_i \leq  k$ set
 $\omega=\alpha_i \circ \omega_L      :  EW(K_i,k)=EW(L_i, k) \lo K_i$.

If $\dim K_i  >  k$ set
$\omega_{k}=\alpha_i \circ \omega_L :    EW_{k}(K_i,k)  =EW(L_i, k) \lo  K_i$
  and we will construct by induction  resolutions
 $\omega_j : EW_j(K_i, k) \lo K_i$, $k+1\leq j\leq \dim K_i$
  such that   $EW_{j}(K_i, k)$ is a subcomplex of
   $EW_{j+1}(K_i, k)$ and $\omega_{j+1}$ extends $\omega_j$  for every $k\leq j  < \dim K_i$.

    Assume that $\omega_j :   EW_j(K_i, k) \lo K_i$, $k\leq j <\dim K_i$
    is constructed.
    For every simplex $\Delta$ of $K_i$ of $\dim=j+1 $   consider  the subcomplex
    ${\omega_j}^{-1}(\Delta)$ of  $EW_j(K_i, k)$.  Enlarge
    ${\omega_j}^{-1}( \Delta)$    by
  attaching  cells  of $\dim =n+1  $
   in order  to kill     the elements of
   $\tor_{\cal F}\pi_n({\omega_j}^{-1}( \Delta))$
   and attaching    cells
   of $ \dim > n+1  $ in order  to  get
    a subcomplex with trivial homotopy groups in
 $\dim   > n$.
 Let   $EW_{j+1}(K_i, k)$ be $EW_j(K_i, k) $ with all the cells  attached   for
 all  $(j+1)$-dimensional simplexes  $\Delta$ of $K_i$ and let
  $\omega_{j+1}:      EW_{j+1} (K_i, k)\lo K_i$ be an extension of
  $\omega_j$
      sending the interior points of the attached cells to the interior of       the corresponding
    $\Delta$.

      Finally denote  $EW(K_i, k)=EW_j(K_i, k)$ and $\omega=\omega_j : EW_j(K_i, k) \lo K_i$
    for $j=\dim K_i$.       Note that since we attach cells only of $\dim >n$,
    the $n$-skeleton of   $EW(K_i, k)$     coincides   with  the $n$-skeleton of $EW(L_i,k)$.

        Let us show that $EW(K_i, k)$ is suitable for $X$.
        Fix a simplex $\Delta$ of $K_i$.
   Since   $\omega^{-1}(\Delta )$ is contractible if $\dim \Delta \leq k$
assume that $\dim \Delta > k$.
         Denote $T= \alpha_i^{-1}(\Delta)$.
      Note that  it follows from the construction that
       $T$  is $(n-1)$-connected,  $H_n(T)$ is $\cal F$-torsion,
       $\omega^{-1}(\Delta )$  is $(k-1)$-connected,
      $\pi_n(    \omega^{-1}(\Delta ))=
 \pi_n (\omega_L^{-1}(T))/ \tor_{\cal F}\pi_n (\omega_L^{-1}(T))$,
          $\pi_j(\omega^{-1}(\Delta ))=0$  for $j \geq  n+1$
    and $\pi_j(\omega^{-1}(\Delta ))=\pi_j (\omega_L^{-1}(T))$ for $ j \leq n-1$.

  Consider the following cases.
   \\\\
  Case 1.  $E=\z_p$.  By (i) of Proposition \ref{p1},
  $\pi_j( \omega_L^{-1}(T))$, $j \leq n-1$
  and $ \pi_n (\omega_L^{-1}(T))/ \tor_{\cal F}\pi_n (\omega_L^{-1}(T))$ are  $p$-torsion.
  Then       $\pi_j(\omega^{-1}(\Delta ))$  is $p$-torsion for $j \leq n$. Therefore
  $\dim_{\pi_j(\omega^{-1}(\Delta ))}  X  \leq \dim_{\z_p} X \leq k$
  for $j \geq k$ and hence by  Proposition \ref{p2}
 $\edim X \leq     \omega^{-1}(\Delta )$.
   \\\\
     Case 2.  $E=\z_{p^\infty}$.  Then by (ii) of Proposition \ref{p1},
     $\pi_j( \omega_L^{-1}(T))$, $j \leq n-1$
  and $ \pi_n (\omega_L^{-1}(T))/ \tor_{\cal F}\pi_n (\omega_L^{-1}(T))$ are  $p$-torsion
  and     $ \pi_k (\omega_L^{-1}(T))$ is   $p$-divisible.
  Then       $\pi_j(\omega^{-1}(\Delta ))$  is $p$-torsion  for $j \leq n$
   and          $\pi_k(\omega^{-1}(\Delta ))$   is
  $p$-divisible.
   Therefore         by the Bockstein theorem and inequalities
   $\dim_{\pi_k(\omega^{-1}(\Delta ))}  X \leq \dim_{\z_{p^\infty} } X  \leq k$   and
  $\dim_{\pi_j(\omega^{-1}(\Delta ))}  X \leq \dim_{\z_{p^\infty} } X +1 \leq k+1$ for $j \geq k+1$.
 Hence by  Proposition \ref{p2}
 $\edim X \leq     \omega^{-1}(\Delta )$.
    \\\\
     Case 3.  $E=\z_{(p)}$.  Then by (iii) of Proposition \ref{p1},
  $\pi_j( \omega_L^{-1}(T))$, $j \leq n-1$   is $p$-local
  and $ \pi_n (\omega_L^{-1}(T))/ \tor_{\cal F}\pi_n (\omega_L^{-1}(T))$ is  $q$-divisible
  for every $q\in \cal P$, $q \neq p$.
  Then       $\pi_j(\omega^{-1}(\Delta ))$, $j \leq n-1$  is $p$-local and
  $\pi_n(\omega^{-1}(\Delta ))$ is $\cal D$-divisible and $\cal F$-torsion free. Therefore
  $\dim_{\pi_j(\omega^{-1}(\Delta ))}  X  \leq k$ for $j\leq n-1$
  and by Proposition \ref{h3}, $\dim_{\pi_n(\omega^{-1}(\Delta ))}  X  \leq n$.
   Hence by  Proposition \ref{p2}
 $\edim X \leq     \omega^{-1}(\Delta )$.
  \\\\
 Case 4. $E=\q$. This case is similar to the previous one.
 \\\\
          Thus we have shown that $EW(K_i, k)$ is suitable for $X$.
       Now   replacing  $K_{i+1}$ by
         $K_j$ with a sufficiently large $j$ we may   assume
    by (i)  of Proposition  \ref{p3}       that there is a combinatorial lifting
         of $h_{i+1}$ to $h'_{i+1} : K_{i+1} \lo  EW(K_i,k)$.
        By (ii) of  Proposition  \ref{p3} we  replace $h'_{i+1}$ by its  cellular approximation
         preserving
  the property of   $h'_{i+1}$ of being a      combinatorial lifting of $h_{i+1}$.

  Let   $\Delta$  be a simplex of   $K_i$ and let
  $\tau :    (\alpha_i\circ \omega_L)^{-1}(\Delta \lo \omega^{-1} (\Delta)$
  be the inclusion.
  Note that from the construction it follows that
   for the induced  homomorphism
   $\tau_* : \pi_n((\alpha_i\circ \omega_L)^{-1}(\Delta)) \lo \pi_n(\omega^{-1} (\Delta))$,
   $\ker    \tau_*$ is $\cal F$-torsion.
  Using this fact and  the reasoning  described in detail in the proof of Theorem 1.2 of \cite{l1}
 one can  construct      from    $K_{i+1}^{[n+1]}$ a CW-complex $L_{i+1}$   by
 replacing some $(n+1)$-simplexes of $K_{i+1}^{[n+1]}$ by $(n+1)$-cells
 attached to the boundary of the replaced simplexes by  a map of degree $\in S({\cal F}) $
such that
    $h'_{i+1}$ restricted
 to    $K_{i+1}^{[n]}$  extends to a map     $g'_{i+1} :  L_{i+1} \lo  EW(L_i, n)$
 such that    $g'_{i+1}, $ $\alpha_{i+1}$, $h_{i+1}$  and
 $\alpha_i \circ \omega_L$  combinatorially commute
 where $\alpha_{i+1}$  is      a projection of $L_{i+1}$ into $K_{i+1}$ taking
 the new cells to the  original ones such that $\alpha_{i+1}$ is 1-to-1 over
 $K^{[n]}_{i+1}$.

   Now define
 $g_{i+1}=\omega_L \circ g'_{i+1} : L_{i+1} \lo L_i$ and
 finally define a simplicial structure on $L_{i+1}$ for which   $\alpha_{i+1}$
 is a combinatorial map.     It is easy to check that the properties (a) and (b)
 are satisfied.
  Since the triangulation of $L_{i+1}$ can be replaced
 by any of its barycentric subdivisions we  may also  assume that  \\

 (c)
 diam$g_{i+1}^j(\Delta) \leq 1/i$ for every simplex
 $\Delta$ in  $L_{i+1}$  and $j \leq i$   \\
     where
 $g^j_i=g_{j+1} \circ g_{j+2} \circ ...\circ g_i: L_i \lo L_j$.   \\

        Denote $Z=\invlim(L_i,g_i)$ and let $r_i : Z \lo L_i$ be the projections.

        Clearly $\dim Z \leq n +1$.
            For constructing
           $L_{i+1}$  we used an arbitrary map $f : N  \lo K(E,k)$ such that
  $E \in \sigma$,          $\dim_{E} X \leq k$,  $ 2 \leq k \leq n$  and  $N$ is a subcomplex
 of $L_i$ containing $L_i^{[k]}$.
By  a standard reasoning described in detail   in the proof of Theorem 1.6 of \cite{l3}
one can show  that  choosing $E$ and $f$  in an appropriate way
           for each $i$ we can achieve
that $\dim_{E} Z \leq k$  for every integer $2\leq k \leq n$ and
 every $E \in \sigma  $  such that      $\dim_{\z_p}  X \leq k $.
Then by the Bockstein theorem
    $\dim_G Z \leq k$ for every group $G$ such that $  \dim_G X \leq k $, $2 \leq k \leq n$.      \\

   By (iii) of Proposition \ref{p3},        the properties (a) and   (b) imply
    that for every   $x \in X$ and  $z\in Z$
   the following  holds:\\

 (d1)
  $ g_{i+1}( \alpha_{i+1}^{-1}(st( p_{i+1} (x) ) )) \subset
   \alpha_{i}^{-1}(st( p_i (x)) )$  and

   (d2)
   $h_{i+1}(st((\alpha_{i+1}\circ r_{i+1})(z)))  \subset        st((\alpha_{i}\circ r_{i})(z))$.
 \\

   Define a  map $r : Z \lo X$ by $r(z)=\cap \{ p_i^{-1}(    st((\alpha_{i}\circ r_{i})(z) ) ): i=1,2,... \}$.
  Then   (d2) implies
   that $r$ is indeed well-defined and continuous.

  The properties (d1) and (d2) also imply  that  for every $x \in X$  \\

 $r^{-1}(x)=\invlim ( \alpha_i ^{-1}(st(p_{i} (x))), g_i |_{\alpha_i ^{-1}(st(p_{i} (x)))})$\\
 where  the map  $ g_i |_{...} $ is considered as a map
      to $\alpha_{i-1} ^{-1}(st(p_{i-1} (x)))$. \\

Since $r^{-1}(x) $ is not empty for every $x \in X$,
 $r$ is a map onto. Fix $x \in X$ and  let us show that $r^{-1}(x)$  satisfies
the conclusions of the theorem.       Since
   $T=\alpha_i ^{-1}(st(p_{i} (x)))$
 is  $(n-1)$-connected  we obtain that
 $r^{-1}(x)$ is approximately  $(n-1)$-connected
as the inverse   limit of $(n-1)$-connected finite simplicial complexes.

Let a group $G$ be such that $\dim_G X \leq n$.
Note that  $H_n(   T )$ is ${\cal F}$-torsion and $G$ is    ${\cal F}$-torsion free.
Then
by  the universal-coefficient theorem $H^n(T;G)={\rm Hom}(H_n(T),G)=0$.
    Thus     ${\Tilde {\Check H}}{}^k(r^{-1}(x); G) =0$ for
     $ k\leq n$ and
 since $\dim_G Z \leq n$,
    ${\Tilde  {\Check H}}{}^k(r^{-1}(x); G) =0$ for  $ k\geq n+1$.
      Hence
  $r$ is $G$-acyclic and  this completes the proof.
\hfill $\Box$

  \end{section}

Department of Mathematics\\
Ben Gurion University of the Negev\\
P.O.B. 653\\
Be'er Sheva 84105, ISRAEL  \\
e-mail: mlevine@math.bgu.ac.il\\\\
\end{document}